\documentclass[12pt,a4paper,reqno]{amsart} 
\usepackage{amssymb}
\usepackage{latexsym}
\usepackage{amsmath}
\usepackage{mathrsfs}
\usepackage[X2,T1]{fontenc}
\usepackage[applemac]{inputenc}
\usepackage{cite}
\usepackage{calc}                   

\newcommand{\Tr}{\operatorname{Tr}}

\newcommand{\scal}[2]{\langle #1,#2\rangle}
\newcommand{\rr}[1]{\mathbf R^{#1}}

\newcommand{\nm}[2]{\Vert #1\Vert _{#2}}

\newcommand{\op}{\operatorname{Op}}

\newcommand{\ep}{\varepsilon}
\newcommand{\fy}{\varphi}
\newcommand{\cdo}{\, \cdot \, }

\newcommand{\wpr}{{\text{\footnotesize $\#$}}}

\newcommand{\eabs}[1]{\langle #1\rangle}     

\newcommand{\vrum}{\vspace{0.1cm}}

\newcommand{\nn}[1]{{\mathbf N}^{#1}}

\newcommand{\maclS}{\mathcal S}

\newcommand{\mascD}{\mathscr D}

\newcommand{\mascF}{\mathscr F}

\newcommand{\mascS}{\mathscr S}
\newcommand{\bsySig}{\boldsymbol \Sigma}

\newcommand{\bsycalS}{\boldsymbol {\mathcal S}}

\newcommand{\WF}{\operatorname{WF}}

\numberwithin{equation}{section}          

\newtheorem{thm}{Theorem}
\numberwithin{thm}{section}

\newcommand{\rubrik}{}
\newtheorem{prop}[thm]{Proposition}
\newtheorem{cor}[thm]{Corollary}
\newtheorem{lemma}[thm]{Lemma}

\theoremstyle{definition}

\newtheorem{defn}[thm]{Definition}

\newcommand{\rubrikdef}{}

\theoremstyle{remark}

\newtheorem{rem}[thm]{Remark}              

\begin{document}

\title{Strong ultra-regularity properties 
for positive elements in the 
twisted convolutions}

\author{Yuanyuan Chen}

\address{Department of Computer science, Physics and Mathematics,
Linn{\ae}us University, Sweden}

\email{yuanyuan.chen@lnu.se}

%
%



\begin{abstract}
We show that positive elements with respect
to the twisted convolutions, belonging to some ultra-test function space of certain 
order at origin, belong to the ultra-test function space of the same order everywhere.
We apply the result to positive semi-definite
Weyl operators.
\end{abstract}

\maketitle

\section{Introduction} \label{sec0}

\par

Several issues in operator theory can be studied by means of
the twisted convolution. For example, composition and positivity
questions can be carried over to related questions
for the twisted convolution product by simple manipulations.
We notice the simple structure of the twisted convolution, since
it essentially consists of a convolution product, disturbed by a
(symplectic) Fourier kernel. It is also common that boundedness
and regularity conditions on operator kernels often
correspond to convenient conditions on related elements in the twisted
convolution. For example, operator kernels which belong to the Schwartz space
$\mascS$, or the Gelfand-Shilov spaces $\maclS _s$ or $\Sigma _s$ of Roumieu
and Beurling types, respectively, carry over to elements in the same class
in the twisted convolution. (See Section \ref{sec1} for notations.)

\par

In \cite{TJ} it is shown that various kinds of singularities for positive elements
with respect to the twisted convolution are attained at the origin.
Furthermore, it is proved that regularity at origin for such elements
impose global regularity and bounedness for these elements and their
Fourier transforms.

\par

More precisely, if $a\in \mascD '$ is positive semi-definite with respect to the twisted
convolution, then it is proved that the following is true:
\begin{enumerate}
\item $a\in \mascS '$ (cf. \cite[Theorem 2.6]{TJ});

\vrum

\item if $\WF _*(a)$ is any wave-front set of $a$ and $(0,Y)\notin \WF _*(a)$,
then $(X,Y)\notin \WF _*(a)$ and $(X,Y)\notin \WF _*(\mascF _\sigma a)$. Here
$\mascF _\sigma$ is the symplectic Fourier transform (cf. \cite[Theorem 4.14]{TJ} and
\cite[Theorem 4.1]{TJ0});

\vrum

\item if $a$ is continuous at origin, then $a$ and its Fourier transform $\widehat a$
are continuous everywhere and belong to $L^2$ (cf. \cite[Theorem 3.13]{TJ});

\vrum

\item if $a\in C^\infty$ near origin, then $a\in \mascS$ (cf. \cite[Theorem 3.13]{TJ});

\vrum

\item if $s\ge 0$, $a\in C^\infty$ near origin and 
\begin{equation}\label{Eq:SmoothGSCond}
|\partial ^\alpha a(0)|\lesssim h^{|\alpha |}\alpha !^s
\end{equation}
for some $h>0$ (for every $h>0$), then $a\in \maclS _s$
($a\in \Sigma _s$) (cf. \cite[Theorem 4.1]{CT}).
\end{enumerate}

\par

We note that if \eqref{Eq:SmoothGSCond} holds true with
$s<1/2$ in (5), then $a$ is trivially equal to $0$, since
the Gelfand-Shilov spaces $\maclS _s$ and $\Sigma _s$
are trivial for such choices of $s$.

\par

In this paper we investigate related questions in background of Pilipovi{\'c}
spaces, $\bsycalS _s$ and $\bsySig _s$ of Roumieu and Beurling type respectively,
a family of function spaces which agrees with corresponding
Gelfand-Shilov spaces when these are non-trivial (cf. \cite{Pil1, Pil2}). We introduce
the so-called twisted Pilipovi{\'c} spaces $\bsycalS _{\sigma ,s}$ and
$\bsySig _{\sigma ,s}$ which are symplectic analogies of Pilipovi{\'c}
spaces, and show that they are homeomorphic to $\bsycalS _s$ and $\bsySig _s$,
respectively. We also show that 
$$
\bsycalS _{\sigma ,s} = \bsycalS _{s} =\maclS _s
$$
when the right-hand side is non-trivial, and similarly for
corresponding spaces of Beurling types.

\par

We consider norm conditions of powers of a second order
partial differential operator $H_\sigma$ and its conjugate. These
operators are symplectic analogies to certain
partial harmonic oscillators. We show that $H_\sigma$ and $\bar
H_\sigma$ commute and can be used to characterize
$\bsycalS _{\sigma ,s}$ and $\bsySig _{\sigma ,s}$ as
\begin{equation}\label{Eq:TwistPilSpChar}
a\in \bsycalS _{\sigma ,s} \ (a\in \bsySig _{\sigma ,s})
\quad
\Leftrightarrow
\quad
\nm {H_\sigma ^{N}\bar H_\sigma^{N}a}{L^\infty}
\lesssim h^{N}(N!)^{4s}
\end{equation}
for some $h>0$ (for every $h>0$). In Section \ref{sec3} we show that
if $a$ is positive semi-definite with respect to the twisted convolution,
then the relaxed condition
$$
|H_\sigma ^{N}\bar H_\sigma ^{N}a(0)|\lesssim h^{N}(N!)^{4s}
$$ 
of the right-hand of \eqref{Eq:TwistPilSpChar} is enough to ensure that
$a$ should belong to $\bsycalS _{\sigma ,s}$ or $\bsySig _{\sigma ,s}$.

\par

\section{Preliminaries}\label{sec1}

\par

In the first part we recall definitions of twisted convolution,
the Weyl quantization and positivity in operator theory, and discuss
basic properties. The verifications are
in general omitted since they can be found in e.{\,}g. \cite{TJ}. Thereafter we
recall the definitions of Gelfand-Shilov and Pilipovi{\'c} spaces and discuss
some properties. Here we also consider related
symplectic analogies of such spaces, defined in terms of Wigner distributions
of Hermite functions, considered by Wong in \cite{Wong, Wo}. Finally we recall some
results in \cite{} on positivity with respect to the twisted convolution.

\par

\subsection{Operators and positivity}\label{subsec1.1}

\par

Let $a$ and $b$ belong to $\mascS (\rr {2d})$, the set of Schwartz functions on $\rr {2d}$.
Then the \emph{twisted convolution} of $a$ and $b$ is given by
$$
(a*_\sigma b)(X) = (2/\pi )^{d/2}\int _{\rr {2d}}a(X-Y)b(Y)e^{2i\sigma (X,Y)}\, dY.
$$
Here $\sigma$ is the symplectic form on $\rr d\times \rr d\simeq \rr {2d}$, given by
$$
\sigma (X,Y) \equiv \scal y\xi -\scal x\eta ,\qquad X=(x,\xi )\in \rr {2d},\ 
Y=(y,\eta )\in \rr {2d}.
$$
The definition of $*_\sigma$ extends in different ways. For example, the map
$(a,b)\mapsto a*_\sigma b$ from $C_0^\infty (\rr {2d})\times C_0^\infty (\rr {2d})$
to $C_0^\infty (\rr {2d})$ is uniquely extendable to a continuous map from
$\mascS '(\rr {2d})\times \mascS (\rr {2d})$ to $\mascS '(\rr {2d})$,
and from $\mascD '(\rr {2d})\times C_0^\infty (\rr {2d})$ to $\mascD '(\rr {2d})$.

\par

There are strong links between the twisted convolution,
and continuity and composition properties in operator theory. This also include
analogous questions in the theory of pseudo-differential operators.

\par

In fact, by straight-forward computations it follows that
\begin{equation}\label{Eq:TwistComp}
A(a*_\sigma b) = (Aa)\circ (Ab),
\end{equation}
where $A$ is the operator defined by the formula
\begin{equation} \label{operatorwithkerAa}
(Aa)(x,y) = (2\pi )^{-d/2}\int _{\rr d}a((y-x)/2,\xi )e^{-i\scal {x+y}\xi }\, d\xi .
\end{equation}
(Here and in what follows we identify operators with their kernels.) We
note that
$$
(Aa)(x,y) = (\mascF ^{-1}(a ((y-x)/2,\cdo )))(-(x+y)),
$$
where $\mascF$ is the Fourier transform on $\mascS '(\rr d)$ which takes the form
$$
\mascF f(\xi ) = \widehat f(\xi )\equiv (2\pi )^{-d/2}\int _{\rr d}f(x)e^{-i\scal x\xi }\, dx
$$
when $f\in \mascS (\rr d)$. Alternatively we may reformulate this identity as
$$
(Aa)(x,y) = (\mascF ^{-1}_2a) ((y-x)/2,-(x+y)),
$$
where $\mascF _2\Phi $ is the partial Fourier transform of $\Phi (x,y)$
with respect to the $y$-variable. Evidently, the mappings $\mascF _2$ and
the pullback which takes $\Phi (x,y)$ into
$$
\Phi ((y-x)/2,-(x+y))
$$
are homeomorphisms on $\mascS (\rr {2d})$ and on $\mascS '(\rr {2d})$,
and unitary on $L^2(\rr {2d})$. Hence similar facts hold true for $A$.

\par

From these mapping properties it follows that if $a\in \mascS '(\rr {2d})$,
then $Aa$ is a linear and continuous operator from $\mascS (\rr d)$
to $\mascS '(\rr d)$. Furthermore, by the kernel theorem of Schwartz
it follows that any linear and continuous operator from $\mascS (\rr d)$
to $\mascS '(\rr d)$ is given by $Aa$, for a uniquely determined $a\in
\mascS '(\rr {2d})$.

\par

At this stage we also note that \eqref{Eq:TwistComp}
remains true, if more generally, $a\in \mascS '(\rr {2d})$ and
$b\in \mascS (\rr {2d})$, which follows by straight-forward computations.

\par

The operator $A$ can also in convenient ways be formulated
in the framework of the Weyl calculus of pseudo-differential operators.
More precisely, the Weyl quantization $\op ^w(a)$ of $a\in \mascS (\rr {2d})$
(the symbol) is the operator from $\mascS (\rr d)$ to $\mascS (\rr d)$
given by
$$
\op ^w(a)f(x) = (2\pi )^{-d}\iint _{\rr {2d}} a((x+y)/2,\xi )
f(y)e^{i\scal {x-y}\xi }\, dyd\xi .
$$
The definition of $\op ^w(a)$ extends in continuous and similar ways as for
$Aa$ to any $\mascS '(\rr {2d})$, and then $\op ^w(a)$ is continuous from
$\mascS (\rr d)$ to $\mascS '(\rr d)$. This extension can also be performed by
the relation
$$
\op ^w(a) = (2\pi )^{-d/2}A(\mascF _\sigma a)
$$
which follows by straight-forward computations.
Here $\mascF _\sigma$ is the symplectic Fourier transform on $\mascS '(\rr {2d})$,
which takes the form
$$
(\mascF _\sigma a)(X) \equiv \pi ^{-d}\int _{\rr {2d}}a(Y)e^{2i\sigma (X,Y)}\, dY
$$
when $a\in \mascS (\rr {2d})$.

\par

From these facts it follow that the Weyl product $\wpr$, defined by
$$
\op ^w(a\wpr b) = \op ^w(a)\circ \op ^w(b)
$$
is given by
$$
a\wpr b =   (2\pi )^{d/2}a*_\sigma (\mascF _\sigma b)
$$
which again links the twisted convolution to compositions in
operator theory.

\medspace

There are also strong links between positivity for the twisted convolution
and positivity in operator theory. We recall that a continuous and linear
operator $T$ from $\mascS (\rr d)$ to $\mascS '(\rr d)$
(from $C^\infty _0(\rr d)$ to $\mascD '(\rr d)$) is called positive
semi-definite, whenever $(Tf,f)\ge 0$ for every $f\in \mascS (\rr d)$
($f\in C^\infty _0(\rr d)$), and then we write $T\ge 0$. Since
$C_0^\infty (\rr d)$ is dense in $\mascS (\rr d)$, it follows
that an operator from $\mascS (\rr d)$ to
$\mascS '(\rr d)$ is positive semi-definite, if it is positive semi-definite
as an operator from $C^\infty _0(\rr d)$ to $\mascD '(\rr d)$.

\par

Positivity for the twisted convolution is defined in an analogous way. That is,
an element $a\in \mascS '(\rr {2d})$ ($a\in \mascD '(\rr {2d})$) is
positive semi-definite with respect to the twisted convolution, whenever
$(a*_\sigma \fy ,\fy )\ge 0$ for every $\fy \in \mascS (\rr {2d})$
($\fy \in C^\infty _0(\rr {2d})$). As above it follows that $a\in \mascS '(\rr {2d})$
is positive semi-definite with respect to $*_\sigma$, if it is positive semi-definite
as an element in $\mascD '(\rr {2d})$.

\par

The following proposition explains the links between positivity in operator theory
and positivity for the twisted convolution. Here $W_{f,g}$ is the Wigner distribution
of $f\in \mascS '(\rr d)$ and $g\in \mascS '(\rr d)$, given by $W_{f,g}\equiv
A^{-1}(\check f\otimes \overline g)$. If $f,g\in \mascS (\rr d)$, then $W_{f,g}$
takes the form
$$
W_{f,g}(x,\xi ) = (2\pi )^{-d/2}\int _{\rr d}f(x-y/2)\overline{g(x+y/2)}e^{i\scal y\xi}\, dy.
$$

\par

\begin{prop}
Let $a\in \mascS '(\rr {2d})$. Then the following conditions are equivalent:
\begin{enumerate}
\item $a$ is positive semi-definite with respect to the twisted convolution;

\vrum

\item $Aa$ is a positive semi-definite operator from $\mascS (\rr d)$ to
$\mascS '(\rr d)$;

\vrum

\item $\op ^w(\mascF _\sigma a)$ is a positive semi-definite
operator from $\mascS (\rr d)$ to $\mascS '(\rr d)$;

\vrum

\item $(\mascF _\sigma a,W_{f,f})\ge 0$ for every $f\in \mascS (\rr d)$.
\end{enumerate}
\end{prop}

\par

\subsection{Gelfand-Shilov spaces}

\par

Let $h,s\in \mathbf R_+$ be fixed. Then $\maclS _{s,h}(\rr d)$
is the set of all $f\in C^\infty (\rr d)$ such that
\begin{equation*}
\nm f{\mathcal S_{s,h}}\equiv \sup \frac {|x^\beta \partial ^\alpha
f(x)|}{h^{|\alpha + \beta |}(\alpha !\, \beta !)^s}
\end{equation*}
is finite. Here the supremum is taken over all $\alpha ,\beta \in
\mathbf N^d$ and $x\in \rr d$.

\par

The set $\mathcal S_{s,h}(\rr d)$ is a Banach space which increases
with $h$ and $s$,
and is contained in $\mascS (\rr d)$.
If $s>1/2$, then $\maclS _{s,h}$ and $\cup _{h>0} \maclS_{1/2, h}$
are dense in $\mathscr S$. Hence, the dual $(\mathcal S_{s,h})'(\rr d)$ of
$\mathcal S_{s,h}(\rr d)$ is a Banach space which contains $\mathscr S'(\rr d)$.

\par

The \emph{Gelfand-Shilov spaces} $\mathcal S_{s}(\rr d)$ and
$\Sigma _s(\rr d)$ are the inductive and projective limits respectively
of $\mathcal S_{s,h}(\rr d)$ with respect to $h>0$. Consequently
\begin{equation*}\label{GSspacecond1}
\mathcal S_s(\rr d) = \bigcup _{h>0}\mathcal S_{s,h}(\rr d)
\quad \text{and}\quad \Sigma _{s}(\rr d) =\bigcap _{h>0}\mathcal
S_{s,h}(\rr d),
\end{equation*}
The space $\Sigma _s(\rr d)$ is a Fr{\'e}chet
space with semi norms $\nm \cdo{\mathcal S_{s,h}}$, $h>0$.
Moreover, $\mathcal S _s(\rr d)\neq \{ 0\}$, if and only if
$s\ge 1/2$, and $\Sigma _s(\rr d)\neq \{ 0\}$,
if and only if $s>1/2$.

\par

If $\ep >0$ and $s>0$, then
$$
\Sigma _s (\rr d)\subseteq \mathcal S_s(\rr d)\subseteq
\Sigma _{s+\ep}(\rr d).
$$

\par

%

\medspace

The \emph{Gelfand-Shilov distribution spaces} $\mathcal S_s'(\rr d)$
and $\Sigma _s'(\rr d)$ are the projective and inductive limits
respectively of $\mathcal S_{s,h}'(\rr d)$.  Hence
$$
\mathcal S_s'(\rr d) = \bigcap _{h>0}\mathcal
S_{s,h}'(\rr d)\quad \text{and}\quad \Sigma _s'(\rr d)
=\bigcup _{h>0} \mathcal S_{s,h}'(\rr d).
$$
By \cite{GS}, $\mathcal S_s'$ and $\Sigma _s'$
are the duals of $\mathcal S_s$ and $\Sigma _s$, 
respectively.

\par

The Gelfand-Shilov spaces and their duals are invariant under
translations, dilations, (partial) Fourier transformations
and under several other important transformations. In fact, by straight-forward
computations it follows that the properties and results in Subsection 1.1
hold true with $\maclS _s$ and $\maclS _s'$ in place of $\mascS$ and
$\mascS '$, respectively, when $s\ge 1/2$, or with $\Sigma _s$ and $\Sigma _s'$
in place of $\mascS$ and $\mascS '$, respectively, when $s>1/2$.

\par

\subsection{The Pilipovi{\'c} spaces}\label{subsec1.2}

\par

We start to consider spaces which are obtained by suitable
estimates of Gelfand-Shilov or Gevrey type when using powers of the harmonic
oscillator $H=|x|^2-\Delta$, $x\in \rr d$.
In general we omit the arguments, since more thorough exposition
is available in e.{\,}g. \cite{TJ4}.

\par

Let $s\ge 0$ and $h>0$. Then $\bsycalS _{\! h,s}(\rr d)$ is 
the Banach space which consists of all
$f\in C^\infty (\rr d)$ such that
\begin{equation}\label{GFHarmCond}
\nm f{\bsycalS _{\! h,s}}\equiv \sup _{N\ge 0}
\frac {\nm{H^Nf}{L^\infty}}{h^N(N!)^{2s}}<\infty .
\end{equation}
If $h_\alpha$
is the Hermite function
\begin{equation} \label{hermitefunt}
h_\alpha (x) = \pi ^{-\frac d4}(-1)^{|\alpha |}
(2^{|\alpha |}\alpha !)^{-\frac 12}e^{\frac {|x|^2}2}
(\partial ^\alpha e^{-|x|^2})
\end{equation}
on $\rr d$ of order $\alpha$, then $Hh_\alpha =(2|\alpha |+d)h_\alpha$.
This implies that $\bsycalS _{\! h,s}(\rr d)$ contains all Hermite functions
when $s>0$, and if
$s=0$ and $\alpha \in \nn d$ satisfies $2|\alpha |+d\le h$,
then $h_\alpha \in \bsycalS _{\! h,s}(\rr d)$.

\par

We let
$$
\bsySig _s(\rr d) \equiv \bigcap _{h>0}\bsycalS _{\! h,s}(\rr d)
\quad \text{and}\quad
\bsycalS _{\! s}(\rr d) \equiv \bigcup _{h>0}\bsycalS _{\! h,s}(\rr d),
$$
and equip these spaces by projective and inductive limit topologies,
respectively, of $\bsycalS _{\! h,s}(\rr d)$, $h>0$. (Cf. \cite{GrPiRo,Pil1,Pil2,TJ4}.)

\par

The space $\bsySig _s(\rr d)${\footnote{The boldface characters  $\bsySig_{s}$, $\bsycalS _{\! s}$, etc. denote
Pilipovi{\'c} spaces, and non-boldface characters 
 $\Sigma_s$, F$\maclS_s$, etc. denote analogous Gelfand-Shilov spaces.}} 
is called the \emph{Pilipovi{\'c} space (of Beurling type) of
order $s\ge 0$ on $\rr d$}.
Similarly, $\bsycalS _{\! s}(\rr d)$ is called the \emph{Pilipovi{\'c} space
(of Roumieu type) of order $s\ge 0$ on $\rr d$}.
Evidently, $\bsySig _0(\rr d)$ is trivially equal to $\{0\}$,
while 
$$
h_{\alpha} \in \bsycalS _{\! s}(\rr d), \quad \text{when} \, s\geq 0
\qquad \text{and} \qquad 
h_{\alpha} \in \bsySig_s(\rr d), \quad \text{when} \, s>0.
$$

\par

The dual spaces of $\bsycalS _{\! h,s}(\rr d)$, $\bsySig _s(\rr d)$
and $\bsycalS _{\! s}(\rr d)$ are denoted by $\bsycalS _{\! h,s}'(\rr d)$,
$\bsySig _s'(\rr d)$ and $\bsycalS _{\! s}'(\rr d)$, respectively.
We have
$$
\bsySig _s'(\rr d) = \bigcup _{h>0} \bsycalS _{\! h,s}'(\rr d)
$$
when $s>0$ and
$$
\bsycalS _{\! s}'(\rr d) = \bigcap _{h>0} \bsycalS _{\! h,s}'(\rr d)
$$
when $s\ge 0$, with inductive respective projective limit topologies
of $\bsycalS _{\! h,s}'(\rr d)$, $h>0$ (cf. \cite{TJ4}).

\par

Let $s>0$ and $\varepsilon >0$. Then
\begin{multline} \label{pilipspinclusion}
\bsycalS _{\! 0}(\rr d) \subseteq \bsySig _{s} (\rr d)\subseteq \bsycalS _{\! s}(\rr d)
\subseteq  \bsySig _{s+\varepsilon}(\rr d) \subseteq \mascS(\rr d)
\\
\subseteq \mascS'(\rr d) \subseteq \bsySig _{s+\varepsilon}'(\rr d)
\subseteq \bsycalS _{\! s}' (\rr d)\subseteq \bsySig _{s}' (\rr d)\subseteq \bsycalS _{\! 0} '(\rr d).
\end{multline}
Furthermore, in \cite{TJ4} it is proved that $\bsycalS_{\! 0}(\rr d)$
consists of all finite linear combinations of Hermite functions,
while $\bsycalS '_{\! 0}(\rr d)$ consists of all formal series 
\begin{equation} \label{pilispacehermiteser}
f=\sum_{\alpha \in \bold N^d}c_{\alpha} h_{\alpha}, 
\qquad
c_{\alpha}=c_{\alpha}(f)=(f, h_{\alpha})_{L^2}.
\end{equation}

\par

The next propositions show that
Pilipovi{\'c} spaces can be characterized by Hermite coefficients
$c_{\alpha}$ given by \eqref{pilispacehermiteser}.
The proofs can be found in \cite{CST, TJ4}. Here $H_1U$ and $H_2U$
are the partial harmonic oscillators given by
\begin{equation} \label{partialharmonicosci}
H_1U(x,y) = (|x|^2-\Delta _x)U(x,y),
\quad 
H_2U(x,y) = (|y|^2-\Delta _y)U(x,y).
\end{equation}

\par

\begin{prop} \label{characGelShi}
Let $s\geq 0$ ($s > 0$) and 
$f \in \bsycalS _{\! 0}'(\rr d)$ be given by
\eqref{pilispacehermiteser}.
Then the following conditions are equivalent:
\begin{enumerate}
\item $f \in \bsycalS _{\! s}(\rr d)$ ($f \in \bsySig_s(\rr d)$);

\vrum

\item $|c_{\alpha}(f)| \lesssim e^{-r|\alpha|^{\frac 1 {2s}}}$
for some $r>0$ (for every $r>0$).
\end{enumerate}
\end{prop}

\par

\begin{prop} \label{characGelShi2}
Let $p,q\in (0,\infty]$, $p_0\in [1,\infty]$, $s\geq 0$ ($s > 0$),
$U \in \bsycalS _{\! 0}'(\rr {2d})$,
and $H_1$ and $H_2$ be given by \eqref{partialharmonicosci}.
Then the following conditions are equivalent:
\begin{enumerate}
\item $U \in \bsycalS _{\! s}(\rr {2d})$ ($U \in \bsySig_s(\rr {2d})$);

\vrum

\item $\nm {H_1^{N_1}H_2^{N_2}U}{L^{p_0}}\lesssim h^{N_1+N_2}(N_1!N_2!)^{2s}$
for some $h>0$ (for every $h>0$);

\vrum

\item $\nm {H_1^{N_1}H_2^{N_2}U}{M^{p,q}}\lesssim h^{N_1+N_2}(N_1!N_2!)^{2s}$
for some $h>0$ (for every $h>0$).
\end{enumerate}
\end{prop}

\par

\begin{rem}
Let $\maclS _s$ and $\Sigma_s$ be the Gelfand-Shilov
spaces of order $s\geq 0$. Then it is proved in \cite{Pil1, Pil2} that
\begin{alignat*}{3}
\bsycalS _{\! s_1} &= \maclS _{s_1},&\quad
\bsySig _{s_2} &= \Sigma _{s_2},& \qquad s_1 &\ge \frac 12,\ s_2 > \frac 12
\intertext{and}
\bsycalS _{\! s_1} &\neq \maclS _{s_1}= \{ 0\},
&\quad \bsySig _{s_2} &\neq \Sigma _{s_2} =\{ 0\} , & \qquad s_1 &<\frac 12,\ 0<s_2
\le \frac 12.
\end{alignat*}
\end{rem}

\par

\begin{rem}
In \cite{TJ4} it is proved that $\bsycalS _{\! s_1}$ and $\bsySig_{s_2}$
are not invariant under dilations when $s_1 < 1/2$ and $s_2 \leq 1/2$.
\end{rem}

\par

\begin{rem}\label{Rem:AltSpaces}
Let the hypothesis in Proposition \ref{characGelShi2} be fulfilled.
By letting $N_1=N_2=N$ we get
\begin{enumerate}
\item[(2)$'$] $\nm {H_1^{N}H_2^{N}U}{L^{p_0}}\lesssim h^{N}N!^{4s}$
for some $h>0$ (for every $h>0$);

\vrum

\item[(3)$'$] $\nm {H_1^{N}H_2^{N}U}{M^{p,q}}\lesssim h^{N}N!^{4s}$
for some $h>0$ (for every $h>0$).
\end{enumerate}
The same arguments as in \cite{CST, TJ4} imply that these conditions are equivalent,
Furthermore, let $\widetilde \bsycalS _{\! s}(\rr {2d})$ ($\widetilde \bsySig _s(\rr {2d})$)
be the set of all $U\in \bsycalS _{\! 0}'(\rr {2d})$ such that
$$
|c_\alpha (U)|\lesssim e^{-r(\eabs {\alpha _1}\eabs {\alpha _2})^{\frac 1{4s}}},
\qquad
\alpha =(\alpha _1,\alpha _2),
$$
for some $r>0$ (for every $r>0$). Then it follows by similar arguments as in
\cite{CST, TJ4} that
$$
(2)' \quad \Leftrightarrow \quad (3)' \quad \Leftrightarrow \quad
U\in \widetilde \bsycalS _{\! s}(\rr {2d})
\
\big (U\in \widetilde \bsySig _{s}(\rr {2d}) \big ).
$$
We note that $\bsycalS _{\! s} \subseteq \widetilde \bsycalS _{\! s}
\subseteq \bsycalS _{2s}$ with strict inclusions.
\end{rem}

\par

\section{Twisted Pilipovi\'c spaces and their properties} \label{sec2}

\par

In this section we introduce twisted Pilipovi\'c spaces as the counter images of
the operator $A$ on Pilipovi{\'c} spaces, and deduce some basic properties.
We also consider their distribution spaces.

\par

We begin with some definitions.

\par

\begin{defn}
The \emph{Hermite-Wong function} of order
$$
\alpha =(\alpha _1,\alpha _2)\in \nn d\times \nn d \simeq \nn {2d}
$$
on $\rr {2d}$ is given by
$$
\varrho _\alpha \equiv A^{-1}(h_{\alpha _1}\otimes h_{\alpha _2})
= A^{-1}(h_{\alpha _1}\otimes \overline {h_{\alpha _2}})
= (-1)^{|\alpha_1|}W_{h_{\alpha_1}, h_{\alpha_2}}.
$$
\end{defn}

\par

The Hermite-Wong functions were studied in different ways by M. W. Wong
in \cite{Wong, Wo}. By the definition it follows that
$$
\varrho _\alpha (X)
=
(2\pi)^{-d/2} \int_{\rr d} h_{\alpha_1}\Big(\frac y 2 - x\Big)
\overline{h_{\alpha_2}\Big(\frac y 2 + x\Big)} e^{i \scal y {\xi}}\, dy,
$$
when $\alpha=(\alpha_1, \alpha_2) \in \nn {2d}$ and $X=(x,\xi )\in \rr {2d}$.

\par

We observe that the Hermite-Wong functions are eigenfunctions
to $\mascF_{\sigma}$. More precisely, we have 
$$
\mascF_{\sigma}\varrho_{\alpha_1, \alpha_2}  
= (-1)^{|\alpha_1|} \varrho_{\alpha_1, \alpha_2},
$$
which follows from the fact that $\mascF_{\sigma}(W_{f,g}) = W_{\check f, g}$
(see e.{\,}g. \cite{Fo, To}. Here $\check f (x) = f(-x)$.

\par

\begin{defn} \label{def:twistpilispace}
Let $s>0$.
\begin{enumerate}
\item The set $\bsycalS '_{\! \sigma, 0}(\rr {2d})$ consists of all formal expansions
\begin{equation} \label{expansiona}
a= \sum_{\alpha}c_{\alpha}\varrho_{\alpha},
\end{equation}
where $\{ c_{\alpha}\} _{\alpha \in \nn {2d}} \subseteq \bold C$.

\vrum

\item The set $\bsycalS _{\! \sigma, 0}(\rr {2d})$ consists of all 
expansions in \eqref{expansiona} such that $c_{\alpha}$
are non-zero for at most finite numbers of $\alpha$.

\vrum

\item The set $\bsycalS _{\! \sigma, s}(\rr {2d})$ ($\bsySig _{\sigma, s}(\rr {2d})$)
consists of all expansions in \eqref{expansiona} such that 
$$
|c_{\alpha}| \lesssim e^{-c|\alpha|^{\frac 1{2s}}}
$$
for some $c>0$ (for every $c>0$).

\vrum

\item The set $\bsycalS '_{\! \sigma, s}(\rr {2d})$ ($\bsySig '_{\sigma, s}(\rr {2d})$)
consists of all expansions in \eqref{expansiona} such that 
$$
|c_{\alpha}| \lesssim e^{c|\alpha|^{\frac 1{2s}}}
$$
for every $c>0$ (for some $c>0$).
\end{enumerate}
\end{defn}

\par

The spaces in Definition \ref{def:twistpilispace} are equipped by topologies
in similar way as for the Pilipov\'c spaces in \cite{TJ4}.

\par

The set $\bsycalS _{\! \sigma, s}(\rr {2d})$ ($\bsySig _{\sigma, s}(\rr {2d})$)
is called the \emph{twisted Pilipovi\'c space of Roumieu type (Beurling type)} 
of order $s$. 
It follows that 
the sets $\bsycalS '_{\! \sigma, s}(\rr {2d})$  
and $\bsySig '_{\sigma, s}(\rr {2d})$
are corresponding distribution spaces,
since similar facts hold true for Pilipovi\'c space \cite{TJ4}.

\par

We extend the definition of $A$ on $\mascS$ by letting
$$
Aa = \sum_{\alpha} c_{\alpha}h_{\alpha}
$$
when $a \in \bsycalS '_{\! \sigma, 0}(\rr {2d})$ is giving by 
\eqref{expansiona}.
It follows that $A$ is a homeomorphism from 
$\bsycalS_{\! \sigma, s}(\rr {2d})$ to $\bsycalS_{\! s}(\rr {2d})$,
from $\bsySig_{\sigma, s}(\rr {2d})$ to $\bsySig_{s}(\rr {2d})$,
and similarly for their duals.
Since it is clear that $A$ is homeomorphism on any Fourier
invariant Gelfand-Shilov spaces, we get
\begin{align*}
\bsycalS _{\! \sigma, s}(\rr {2d}) = \bsycalS _{\!  s}(\rr {2d})
= \maclS_s(\rr {2d}), \quad \text{when} \, s\geq 1/2
\intertext{and}
\bsySig _{\sigma, s}(\rr {2d}) = \bsySig _{s}(\rr {2d})
= \Sigma_s(\rr {2d}), \quad \text{when} \, s> 1/2,
\end{align*}
and similarly for corresponding distribution spaces.

\par

\begin{rem}
Let $a \in \bsycalS '_{\! \sigma, 0}(\rr {2d})$ be as 
in \eqref{expansiona}. Since $A$ is a homeomorphism
on $\mascS(\rr {2d})$ and on $\mascS'(\rr {2d})$,
it follows from \cite{TJ1} that $a$ belongs to $\mascS(\rr {2d})$ 
if and only if $c_{\alpha} \lesssim \eabs x ^{-N}$
for every $N \geq 0$. In the same way, 
$a \in \mascS'(\rr {2d})$ if and only if 
$c_{\alpha} \lesssim \eabs x ^{N}$
for some $N \geq 0$.
\end{rem}

\par

Next we discuss the partial harmonic oscillators $H_1$
and $H_2$ in Proposition \ref{characGelShi2}, 
and their counter images under
the operator $A$.
We let $H_{\sigma}$ be the operator on $\mascS(\rr {2d})$,
given by 
$$
H _{\sigma} = (|X|^2 -\frac 1 4\Delta _X) +\scal {\xi}{D_x} -\scal{x}{D_{\xi}},
\quad X=(x, \xi) \in \rr {2d},
$$
and we let $T_{\sigma}=H _{\sigma} \circ \bar H _{\sigma}$.
Here we note that
$$
 \bar H _{\sigma}
=  (|X|^2 -\frac 1 4\Delta _X) -\scal {\xi}{D_x} +\scal{x}{D_{\xi}}.
$$

\par

The following lemma explains some spectral properties
of the considered operators.

\par

\begin{lemma} \label{symharmonicopeigen}
Let $s\geq 0$.
Then the following is true:
\begin{enumerate}
\item the Hermite-Wong functions $\varrho_{\alpha}$ are  
eigenfuctions to $H_{\sigma}$, $\bar H_{\sigma}$
and $T_{\sigma}$, and 
\begin{align} \label{eigenfunctoTsigma}
H_{\sigma} \varrho_{\alpha_1, \alpha_2} 
= (2|\alpha_1|+d) \varrho_{\alpha_1, \alpha_2}, \quad 
\bar H_{\sigma} \varrho_{\alpha_1, \alpha_2} 
= (2|\alpha_2|+d) \varrho_{\alpha_1, \alpha_2}, 
\intertext{and}
T_{\sigma} \varrho_{\alpha_1, \alpha_2} 
= (2|\alpha_1|+d) (2|\alpha_2|+d) \varrho_{\alpha_1, \alpha_2}; \nonumber
\end{align}

\vrum

\item $H_{\sigma}$ and $\bar H_{\sigma}$ 
restrict to homeomorphisms on 
$\bsycalS _{\! \sigma, s}(\rr {2d})$ 
and on $\bsySig _{\sigma, s}(\rr {2d})$;
\vrum

\item the definitions of 
$H_{\sigma}$ and $\bar H_{\sigma}$ 
extend uniquely to homeomorphisms on 
$\mascS'(\rr {2d})$,
$\bsycalS '_{\! \sigma, s}(\rr {2d})$ 
and on $\bsySig '_{\sigma, s}(\rr {2d}) $.
\end{enumerate}
\end{lemma}

\par

For the proof, we shall make use of the operators
\begin{equation*} 
\begin{alignedat}{2}
Z_{1,j} &= \frac 1 {2} \partial_{z_j} + \overline z_j, & \qquad
\widetilde Z_{1,j} &= \frac 1 {2} \partial _{\overline z_j} - z_j,
\\[1ex]
Z_{2,j} &=  \frac 1 {2} \partial _{\overline z_j} + z_j, & \qquad
\widetilde Z_{2,j} &=\frac 1 {2} \partial_{z_j} - \overline z_j,
\end{alignedat}
\end{equation*}
where
\begin{equation*} 
\begin{alignedat}{2}
z_j &= x_j+ i\xi_j, & \qquad
\overline z_j &= x_j - i\xi_j,
\\[1ex]
\partial_{z_j} &=\partial_{x_j} - i\partial_{\xi_j}, & \qquad
\partial _{\overline z_j} &= \partial_{x_j} + i\partial_{\xi_j},
\end{alignedat}
\end{equation*}
(see \cite[Section 22]{Wo}).
By similar arguments as in the proof of Theorem 22.1 in \cite{Wo}
we get
\begin{equation} \label{anniandcreatop}
\begin{alignedat}{1} 
Z_{1,j}\varrho_{\alpha_1, \alpha_2} & = (2|\alpha_{2,j}|)^{1/2}\varrho_{\alpha_1, \alpha_2-e_j}, 
\\
\widetilde Z_{1,j}\varrho_{\alpha_1, \alpha_2} & = -(2|\alpha_{2,j}|+2)^{1/2} \varrho_{\alpha_1, \alpha_2+e_j},
\\
Z_{2,j}\varrho_{\alpha_1, \alpha_2} & = -(2|\alpha_{1,j}|)^{1/2}\varrho_{\alpha_1-e_j, \alpha_2}, 
\\
\widetilde Z_{2,j}\varrho_{\alpha_1, \alpha_2} & =  (2|\alpha_{1,j}|+2)^{1/2}\varrho_{\alpha_1+e_j, \alpha_2},
\end{alignedat}
\end{equation}
where $e_1, \dots, e_d$ is the standard basis in $\rr d$, i.e, 
$e_j=(\delta_{1,j}, \dots, \delta_{d,j})$, $j=1, \dots, d$, 
and $\delta_{i,j}$ is the Kroniker's delta function.

\par

In view of \eqref{anniandcreatop}, the operators $Z_{1,j}$ and $Z_{2,j}$
can be considered as symplectic analogies of annihilation operators,
$\widetilde Z_{1,j}$ and $\widetilde Z_{2,j}$ as symplectic analogies
of creation operators.

\par

\begin{proof}
First we prove (1). 
By straight-forward computations, we obtain
\begin{align*}
H_{\sigma}= -\frac 1 2 \Big( \sum_j Z_{2,j}\widetilde Z_{2,j} + \widetilde Z_{2,j}Z_{2,j}\Big)
\intertext{and}
\bar H_{\sigma} = -\frac 1 2 \Big(\sum_j Z_{1,j}\widetilde Z_{1.j} + \widetilde Z_{1,j}Z_{1,j}\Big).
\end{align*}
Hence, by \eqref{anniandcreatop} we get
\begin{align*}
H_{\sigma} \varrho_{\alpha_1, \alpha_2}= (2|\alpha_1|+d) \varrho_{\alpha_1, \alpha_2},
\intertext{and}
\bar H_{\sigma} \varrho_{\alpha_1, \alpha_2}= (2|\alpha_2|+d) \varrho_{\alpha_1, \alpha_2},
\end{align*}
and (1) follows.

\par

By \eqref{eigenfunctoTsigma}, 
it follows that $H_{\sigma}$ and $\bar H_{\sigma}$ 
restrict to homeomorphisms on 
$\bsycalS _{\! \sigma, s}(\rr {2d})$ 
and on $\bsySig _{\sigma, s}(\rr {2d})$,
which gives (2).

\par

If $a \in \bsycalS '_{\! \sigma, s}(\rr {2d})$
and $b \in \bsycalS _{\! \sigma, s}(\rr {2d})$.
We now let $H_{\sigma}$ be defined by
$$
(H_{\sigma}a, b)_{L^2} = (a, \bar H_{\sigma}b)_{L^2},
$$
as usual, which extends the definitions of 
$H_{\sigma}$ and $\bar H_{\sigma}$ to
$\bsycalS '_{\! \sigma, s}(\rr {2d})$.
The extensions of these operators to 
$\bsySig '_{\sigma, s}(\rr {2d})$ and  
$\mascS'(\rr {2d})$ are performed 
in similar ways.
By \eqref{eigenfunctoTsigma}, 
it follows that these extensions are unique.
\end{proof}

\par

The next lemma shows important links between the latter 
operators and partial harmonic oscillators. 

\par

\begin{lemma} \label{decomofATNa}
Let $H_1$ and $H_2$ be as in Proposition \ref{characGelShi2},
and let $a \in \bsycalS _{\! \sigma, s}(\rr {2d})$.
Then $H_{\sigma}$ and $\bar H_{\sigma}$
are commuting to each other, and  
\begin{equation} \label{relationAHT}
A(H_{\sigma}^{N_1}\bar H_{\sigma}^{N_2}a) 
= H_1^{N_1} H_2^{N_2}(Aa),
\end{equation}
for every interger $N_1, N_2 \geq 0$.
In particular, if $\{ f_k\} _{k=1} ^{\infty}$ and $\{ g_k\} _{k=1} ^{\infty}$
are sequences in $l^2(\bold N; L^2(\rr d))$, and $a$ is given
by 
$$
a=\sum_{k=0} ^{\infty} A^{-1}(f_k \otimes \overline{g_k}),
$$
then 
$$
A(T_{\sigma}^Na) = \sum_{k=0}^{\infty}(H^Nf_k) \otimes (\overline{H^Ng_k}),
$$
where the series convergences in $\mascS'(\rr {2d})$.
\end{lemma}

\par

\begin{proof}
The commutation between $H_{\sigma}$ and $\bar H_{\sigma}$
follows if we prove \eqref{relationAHT}.
We recall the 
operators 
\begin{equation} 
\begin{alignedat}{2}
P_j &= \frac 1 {2i} \partial_{\xi _j} - x_j, & \qquad
\Pi _j &= \frac 1 {2i} \partial _{x_j} + \xi _j,
\\[1ex]
T _j &= \frac 1 {2i} \partial _{\xi _j} + x_j, & \qquad
\Theta _j &= \frac 1 {2i} \partial _{x_j} -\xi _j,
\end{alignedat}
\end{equation}
and the relations
\begin{equation}\label{eq.xxMotiv}
\begin{alignedat}{2}
A(P_j ^2a) &= x_j ^2 Aa,&\qquad A(\Pi _j ^2a) &= -\partial _{x_j} ^2(Aa),
\\[1ex]
A(T_j ^2 a) &= y_j ^2 Aa, &\qquad  
A(\Theta _j ^2a) &=-\partial_{y_j} ^2 (Aa),
\end{alignedat}
\end{equation}
from \cite[Theorem 4.1]{CT}.

\par

By straight-forward computations we get
$$
(x_j^2-\partial_{x_j}^2)(Aa) 
= A((P _j^2 + \Pi _j^2)a)
=A(H_{\sigma, j}a), 
$$
where $H_{\sigma, j}= (X _j^2 -\frac 1 4\Delta _{X_j}) + \xi_j D_{x_j}-x_jD_{\xi_j}$.

\par

Summing up over all $j$ gives
$$
H_1(Aa) = A(H_{\sigma}a).
$$
In the same way we get 
$$
H_2(Aa) = A(\bar H_{\sigma}a),
$$
and the result follows by induction.
\end{proof}

\par

From these mapping properties,
Proposition \ref{characGelShi2}
can now be carried over to the case of 
twisted Pilipovi\'c spaces as follows.

\par

\begin{prop}
Let $p, q \in (0, \infty]$ and $p_0 \in [1, \infty]$
and let $s\geq 0$ ($s>0$).
Then the following conditions are equivalent.
\begin{enumerate}
\item $a \in \bsycalS _{\! \sigma, s}(\rr {2d})$ ($a \in \bsySig _{\sigma, s}(\rr {2d})$);

\vrum

\item $\nm {H_{\sigma}^{N_1} \bar H_{\sigma}^{N_2}a}{L^{p_0}}\lesssim h^{N_1+N_2}(N_1!N_2!)^{2s}$
for some $h>0$ (for every $h>0$);

\vrum

\item $\nm {H_{\sigma}^{N_1} \bar H_{\sigma}^{N_2}a}{M^{p, q}}\lesssim h^{N_1+N_2}(N_1!N_2!)^{2s}$
for some $h>0$ (for every $h>0$).
\end{enumerate}
\end{prop}

\par

\begin{proof}
Let $U=Aa$. Since 
$M^{p_1}(\rr {2d}) \subseteq M^{p, q}(\rr {2d}) \subseteq M^{p_2}(\rr {2d})$,
when $p_1 = \min (p, q)$ and $p_2 = \max (p, q)$, 
we may assume that $p=q$.

\par

Since $A$ is a homeomorphism on $M^p(\rr {2d})$,
we get
$$
\nm {H_{\sigma}^{N_1} \bar H_{\sigma}^{N_2}a}{M^{p}} 
= \nm {A(H_{\sigma}^{N_1} \bar H_{\sigma}^{N_2}a)}{M^{p}}
= \nm {H_1^{N_1} H_2^{N_2}U}{M^{p}},
$$
and the equivalence between (3) and Proposition \ref{characGelShi2}
(3) follows. The equivalence between (1) and (3) now follows from 
Proposition \ref{characGelShi2} and the fact that $A$
is a homeomorphism from $\bsycalS _{\! \sigma, s}(\rr {2d})$ 
to $\bsycalS _{\! s}(\rr {2d}) $. 

\par

Finally by the embeddings
$$
M^1(\rr {2d}) \subseteq L^{p_0}(\rr {2d}) \subseteq M^{\infty}(\rr {2d}),
$$
the equivalence between (2) and (3) now follows.
\end{proof}

\par

\begin{cor} \label{aimplyTsigma}
If $s \geq 0$ and $a \in \bsycalS _{\! \sigma, s}(\rr {2d})$
($a \in \bsySig _{ \sigma,s}(\rr {2d})$), then
\begin{equation} \label{opTsigmaintwistedps}
\nm {T_{\sigma}^Na} {L^{\infty}} \lesssim h^{2N}(N!)^{4s},
\end{equation}
holds for some $h>0$ (for every $h>0$).
\end{cor}

\par

Remark \ref{Rem:AltSpaces} and Lemma \ref{decomofATNa} show that
\eqref{opTsigmaintwistedps} is necessary but not sufficient in order for
$a \in \bsycalS _{\! \sigma, s}(\rr {2d})$ or $a \in \bsySig _{ \sigma,s}(\rr {2d})$.

\par

\section{Twisted Pilipovi{\'c} space property for positive
elements with respect to  the twisted convolution} \label{sec3}

\par

We study positive elements with respect to 
twisted convolution in $\mascS'$,
having the twisted Pilipovi{\'c} space regularities near the origin.
We show that such elements are in 
$\bsycalS _{\! \sigma, s}$ or in $\bsySig _{\sigma, s}$.

\par

The following theorem shows that the condition
of the form \eqref{opTsigmaintwistedps}
at origin is sufficient that the converse
of Corollary \ref{aimplyTsigma} holds when dealing with positive 
semi-definite elements with respect to the
twisted convolution.

\par

\begin{thm} \label{sympilispace}
Let $s\geq 0$, $a\in \mascS '(\rr {2d})$ and 
$(a \ast _{\sigma}\psi, \psi) \geq 0$ for every $\psi \in \mascS (\rr {2d})$. 
If
$$
(T_{\sigma}^Na)(0,0) \lesssim h^{2N}(N!)^{4s},
$$ 
holds for some $h>0$ (for every $h>0$), 
then $a \in \bsycalS _{\! \sigma,s}(\rr {2d})$
($a \in \bsySig _{ \sigma,s}(\rr {2d})$).
\end{thm}

\par

\begin{proof}
By the assumption, we may write $a=\sum_k A^{-1}(f_k \otimes \overline{f_k})$.
By Lemma \ref{decomofATNa}, we obtain 
$$
A(T_{\sigma}^Na) 
=
\sum_k (H^Nf_k \otimes \overline{H^Nf_k}),
$$
for some sequence $\{ f_k\}_{k=0}^{\infty}$.

Let $K=\sum_kf_k \otimes \overline{f_k}$ be the kernel of $Aa$.
Then
\begin{multline*}
\|H_1^N H_2^N K \|_{L^2} 
\leq 
\| H_1^N H_2^N K \|_{\Tr} 
\\
=
\| A(T_{\sigma}^Na) \|_{\Tr} 
=
\sum_k \|H^Nf_k \|_{L^2} ^2=\big( {\pi} / 2 \big)^{d/2}(T_{\sigma}^Na)(0,0).
\end{multline*}
Thus by the assumption, we get
$$
\|H_1^N H_2^N K \|_{L^2}  \lesssim h^{2N}(N!)^{4s},
$$
for some $h>0$ (for every $h>0$),
giving that $K \in \bsycalS _{\! s}(\rr {2d})$ ($K \in \bsySig _{s}(\rr {2d})$)
in view of Proposition \ref{characGelShi2} and Remark \ref{Rem:AltSpaces}.
Hence $a \in \bsycalS _{\! \sigma,s}(\rr {2d})$
($a \in \bsySig _{ \sigma,s}(\rr {2d})$).

\end{proof}

\par

\begin{prop}
Let $s \geq 0$ be real, 
and let $a \in \mascS '(\rr {2d})$ be such that
$\op ^{\omega}(a) \geq 0$. If
\begin{equation} \label{condsympilispacepsu}
(T ^N_{\sigma} (\mascF _{\sigma}a))(0) \lesssim h ^{2N} (N!) ^{4s},
\end{equation}
holds for some $h>0$ (for every $h>0$), 
then $a \in \bsycalS _{\! \sigma,s}(\rr {2d})$
($a \in \bsySig _{ \sigma,s}(\rr {2d})$).
\end{prop}

\par

\par

\end{document}